\documentclass[twoside,11pt]{article}

\usepackage{amsmath,amsfonts}
\usepackage{indentfirst}
\usepackage{latexsym,bm}

\newcommand{\D}{\displaystyle}

\begin{document}
\title{\bf Asymptotic Behavior of a Viscous Liquid-Gas Model  with Mass-Dependent
Viscosity and Vacuum}
\author{{Q{\sc ingqing} L{\sc iu} \ \ \ \ \ \ \ \ \ \ \ \ \ C{\sc hangjiang} Z{\sc hu}\thanks{Corresponding author.
Email: cjzhu@mail.ccnu.edu.cn}}   \\  \\
{The Hubei Key Laboratory of Mathematical Physics} \\
{School of Mathematics and Statistics }
\\{Huazhong Normal University, Wuhan 430079, P.R. China}}
\date{}

\maketitle

\begin{abstract}

In this paper, we consider two classes of free boundary value
problems of a viscous two-phase liquid-gas model relevant to the
flow in wells and pipelines with mass-dependent viscosity
coefficient. The liquid is treated as an incompressible fluid
whereas the gas is assumed to be polytropic. We obtain the
asymptotic behavior and decay rates of the mass functions $n(x,t)$,\
$m(x,t)$ when the initial masses are assumed to be connected to
vacuum both discontinuously and continuously, which improves the
corresponding result about Navier-Stokes equations in \cite{Zhu}.
\end{abstract}
 \vspace{3mm}

\textbf{Key words:} Two-phase flow model, mass-dependent viscosity,
vacuum, asymptotic behavior.

\vspace{3mm}
 {\textbf{AMS Subject Classification (2000):} 76T10,
35R35, 35B40.}

\section*{Contents}

\noindent 1.  Introduction \dotfill 2

\noindent 2.  Reformulation of the problems and the main results
\dotfill 4

\noindent 3.  A \emph{priori} estimates and the asymptotic behavior
of the mass functions \dotfill 8

3.1. Uniform {\it a priori} estimates \dotfill 8

3.2. Asymptotic behavior of $cQ(m)$ \dotfill 14

\noindent 4. Decay rates  of the mass functions  \dotfill 15

 \noindent Acknowledgement \dotfill 22

 \noindent References \dotfill 23

\section{Introduction }

Consider the one dimensional liquid-gas two-phase model with
viscosity which can be written in Eulerian coordinates as (cf.
\cite{Evje-Karlsen,Evje-Flatten-Friis}):
\begin{eqnarray}\label{1.1}
\left \{
\begin{array}{l}
\partial_{t}[\alpha_{g} \rho_{g}]+  \partial_{x}[\alpha_{g}
\rho_{g}u_{g}] = 0 ,\\[3mm]
\partial_{t}[\alpha_{l} \rho_{l}] +  \partial_{x}[\alpha_{l}
\rho_{l}u_{l}] =0,\\[3mm]
 \partial_{t}[\alpha_{g}\rho_{g}u_{g}+  \alpha_{l}\rho_{l}u_{l}]+\partial_{x}[\alpha_{g}\rho_{g}u_{g}^{2}+
 \alpha_{l}\rho_{l}u_{l}^{2}+P]=-q+\partial_{x}[\varepsilon
 \partial_{x}u_{mix}],
\end{array}
\right.
\end{eqnarray}
where $u_{mix}=\alpha_{g}u_{g}+\alpha_{l}u_{l}$ and the unknown
 variables $\alpha_{g},\alpha_{l} \in [0,1]$  denote volume
 fractions satisfying the fundamental relation:
\begin{eqnarray}
 \alpha_{g}+\alpha_{l}=1.
\end{eqnarray}
 Furthermore, the other unknown
 variables $\rho_{g},\ \rho_{l},\ u_{g},\ u_{l}$ denote gas density, liquid
 density, velocities of gas and liquid respectively, whereas $P$ is  the common
 pressure for both phases, $q$ presents external forces, like
 gravity and friction, and $\varepsilon>0$ denotes viscosity.

We focus on a simplified model as in
\cite{Evje-Karlsen,Evje-Flatten-Friis} ob\/tained by assuming that
fluid velocities are equal, i.e., $u_{g}=u_{l}=u$ and  neglecting
the exter\/nal forces, i.e., $q=0$. In addition, we  neglect the gas
phase effects in the mixture momentum conservation equation
$(1.1)_{3}$. This motivation is from the fact that liquid density is
much higher than the gas density, generality speaking,
$\rho_{l}/\rho_{g}=O(10^{3})$. Thus, we can obtain the following
simplified model:
\begin{eqnarray}\label{1.3}
\left \{
\begin{array}{l}
\partial_{t}[\alpha_{g} \rho_{g}]+  \partial_{x}[\alpha_{g}
\rho_{g}u] = 0, \\[3mm]
\partial_{t}[\alpha_{l} \rho_{l}] +  \partial_{x}[\alpha_{l}
\rho_{l}u] =0,\\[3mm]
 \partial_{t}[ \alpha_{l}\rho_{l}u]+\partial_{x}[
 \alpha_{l}\rho_{l}u^{2}+P]=\partial_{x}[\varepsilon
 \partial_{x}u].
\end{array}
\right.
\end{eqnarray}

As in \cite{Evje-Karlsen}, we assume  the liquid is incompressible,
i.e., $\rho_{l}=${\it constant} and the gas is polytropic, i.e.,
$$
P=C\rho_{g}^{\gamma}, \ \ \ \gamma>1, \ \ \ C>0. \eqno(1.4)
$$
Let
$$
n=\alpha_{g}\rho_{g},\ \ \ \ \ m=\alpha_{l}\rho_{l}.\eqno(1.5)
$$
By (1.2), (1.3), (1.4) and (1.5), we have
$$
\left \{
\begin{array}{l}
\partial_{t}n+  \partial_{x}(nu) = 0 ,\\[3mm]
\partial_{t}m+  \partial_{x}(mu) = 0 ,\\[3mm]
 \partial_{t}(mu)+\partial_{x}(mu^{2}+P)=\partial_{x}(\varepsilon
 \partial_{x}u),
\end{array}
\right. \eqno(1.6)
$$
and
$$P(n,m)=C\rho_{l}^{\gamma} \left (\frac{n}{\rho_{l}-m}\right)^{\gamma}
= A \left (\frac{n}{\rho_{l}-m}\right)^{\gamma}, \eqno(1.7)$$ where
$A=C\rho^{\gamma}_l$. Moreover the viscosity coefficient is taken as
the form (cf. \cite{Evje-Karlsen, Evje-Flatten-Friis}):
$$\varepsilon=\varepsilon(n,m)=B \frac{n^{\beta}}{(\rho_{l}-m)^{\beta+1}},\ \ \ \ B>0,\ \beta>0,\eqno(1.8)$$
or
$$\varepsilon=\varepsilon(m)=B \frac{m^{\beta}}{(\rho_{l}-m)^{\beta+1}},\ \ \ \ B>0,\  \beta>0.\eqno(1.9)$$
In the following, without loss of generality, we consider only the
form $(1.8)$.

In this paper, we will consider the following two classes of the
free boundary value problems of the system (1.6):

(1) The initial masses connect to vacuum discontinuously:

     The boundary conditions are given as
     $$
     \left\{
     \begin{array}{l}
     (-P(m,n)+\varepsilon(n,m)\partial_{x}u)(a(t)^{+},t)=0,\\[3mm]
      (-P(m,n)+\varepsilon(n,m)\partial_{x}u)(b(t)^{-},t)=0, \ \ \
      t\geq 0,
     \end{array}
     \right.
     \eqno(1.10)
     $$
and the initial data are given as
$$n(x,0)=n_{0}(x)>0,\ \ m(x,0)=m_{0}(x)>0,\ \ u(x,0)=u_{0}(x),
\ \ x \in [a,b].\eqno(1.11)$$

     (2) The initial  masses connect to vacuum continuously:

     The boundary conditions are given as
     $$n(a(t),t)=n(b(t),t)=0,\ \ \ \ \ m(a(t),t)=m(b(t),t)=0, \ \ \ t\geq 0, \eqno(1.12)$$
and the initial data are given as
$$n(x,0)=n_{0}(x)>0, \ \ \ m(x,0)=m_{0}(x)>0, \ \ \ u(x,0)=u_{0}(x),
\ \ \ x \in (a,b),\eqno(1.13)$$ and
$n_{0}(a)=n_{0}(b)=m_{0}(a)=m_{0}(b)=0$.

Here $-\infty<a<b<\infty$, $a(t)$ and $b(t)$ are the free boundaries
defined by
 $$
\left\{\begin{array}{l} \displaystyle \frac{da(t)}{dt}=u(a(t), t), \
\ \ t>0,
\\ [3mm] a(0)=a,
\end{array}
\right. \eqno(1.14)
$$
and
$$
\left\{\begin{array}{l} \displaystyle \frac{db(t)}{dt}=u(b(t), t), \
\ \ t>0,
\\ [3mm] b(0)=b,
\end{array}
\right. \eqno(1.15)
$$
which are the interfaces separating the gas from the vacuum.

Let's  first  review some of the previous works in this direction.
When the viscosity coefficient $\varepsilon$ was taken as the form
$(1.9)$ and the initial masses connected to vacuum discontinuously,
Evje and Karlsen in \cite{Evje-Karlsen} got the global existence and
uniqueness of weak solutions when $\beta \in (0,\frac{1}{3})$ by
energy method. This result was later generalized to the case when
$\beta\in (0,1]$ by Yao and Zhu in \cite{Yao-Zhu2}. When the
viscosity coefficient $\varepsilon$ was taken as the form $(1.8)$
and the initial masses connected to vacuum continuously, Evje,
Flatten and Friis in \cite{Evje-Flatten-Friis} proved the global
existence of the weak solutions under some approximate assumptions
on $n_{0}(x),\ m_{0}(x),\ u_{0}(x),$ and $\displaystyle
\frac{n_{0}(x)}{m_{0}(x)}=c_{0}(x)$ when $\beta \in
(0,\frac{1}{3})$. This result was  later generalized to the case
$\beta \in (0,1)$ by Yao and Zhu in \cite{Yao-Zhu1} when the initial
masses connected to vacuum continuously and the viscosity
coefficient was a constant.

It is necessary for us  to illustrate that the main methods used to
obtain our results are similar to those  in  \cite{Zhu,Guo-Zhu},
because we have used the variable  transformations as in
\cite{Evje-Karlsen,Evje-Flatten-Friis}, by  which we can rewrite our
problem into (2.7)-(2.12) similar to the model in single-phase
Navier-Stokes equations. In view of this, let's review some of the
relevant works about single-phase Navier-Stokes equations with
density-dependent viscosity and vacuum. When the initial density
connected to vacuum discontinuously, the global existence of weak
solutions for isentropic flow was obtained by  Okada,
Matu$\check{\rm s}\dot{\rm u}$-Ne$\check{\rm c}$asov$\acute{\rm a}$,
Makino in \cite{Okada2} for $\mu(\rho)=\rho^{\theta}$,
$0<\theta<\frac{1}{3}$, by Yang, Yao and Zhu in \cite{Yang-Yao-Zhu}
for $0<\theta<\frac{1}{2}$ and by Jiang, Xin and Zhang in
\cite{Jiang-Xin-Zhang} for $0<\theta<1$. Qin, Yao and Zhao in
\cite{Qin-Yao-Zhao} extended the results in \cite{Okada2,
Yang-Yao-Zhu, Jiang-Xin-Zhang} to the case $0<\theta\leq1$.
Recently, Zhu in \cite{Zhu} investigated the asymptotic behavior and
decay rate estimates about the density function $\rho(x,t)$ by
overcoming some new difficulties which came from the appearance of
the boundary layers. When the initial density connected to vacuum
continuously, the local existence of weak solutions was obtained in
\cite{Yang-Zhao} by Yang and Zhao. The global existence of weak
solutions was given in \cite{Yang-Zhu} by Yang and Zhu for
$0<\theta<\frac{2}{9}$, and later was improved in
\cite{Vong-Yang-Zhu} for $0<\theta<\frac{1}{3}$ by Vong, Yang and
Zhu, in
 \cite{Fang-Zhang1} for $0<\theta<\frac{1}{2}$ and in
 \cite{Fang-Zhang2} for $0<\theta<1$ by Fang and Zhang. Guo and Zhu in  \cite{Guo-Zhu} firstly  studied  the asymptotic
  behavior and the decay rate of the density function $\rho(x,t)$
  with respect to  the time $t$ for any $\theta>0$ based on the
  following new mathematical entropy inequality, which was obtained first by Kanel in \cite{Kanel} for one-dimensional case and
Bresch, Desjardins, Lin and Mellet, Vasseur for multi-dimensional
case, cf. \cite{Bresch1, Bresch2, Mellet}:
  $$\arraycolsep=1.5pt
  \begin{array}{rl}
  &\displaystyle \int^{b(t)}_{a(t)} \left \{  \frac{1}{2}\rho u^{2}+ u(\rho^{\theta})_{x}
  +\frac{1}{2}\rho^{2\theta-3}\rho_{x}^{2}+\frac{\rho^{\gamma-1}}{\gamma-1}\right \}dx \\[3mm]
  +&\displaystyle \int^{t}_{0}\int^{b(t)}_{a(t)}\frac{4\theta\gamma}{(\gamma+\theta)^{2}}
  \rho^{\gamma+\theta-3}\rho_{x}^{2}dxdt\leq C,
  \end{array}$$
where $C $ is a uniform constant independent of $t$.

In this paper, we show that the masses $n$ and $m$ tend to zero as
time goes to infinity. Moreover, we can obtain a stabilization rate
estimates of the mass  functions for any $\beta>0$ as
$t\rightarrow\infty$.

The rest of this paper is organized as follows. In Section 2, we
reformulate the two free boundary value problems (1.6), (1.10),
(1.11) and (1.6), (1.12), (1.13) into the two fixed boundary value
problems by introducing the Lagrangian coordinates and using the
variable transformations. Then we state the main theorems of this
paper. In Section 3, we derive some crucial uniform estimates for
studying the asymptotic behavior and the decay rate estimates about
the mass functions. In Section 4, the decay rate estimates on the
mass functions  will be given by introducing a new function $w(x,t)$
in \cite{Nagasawa} by Nagasawa.

\section{Reformulation of the problems and the main results}

To solve the two free boundary problems above, it is convenient to
convert the free boundaries to the  fixed boundaries by using
Lagrangian coordinates. To do this, let
$$\xi=\int_{a(t)}^{x}m(y,t)dy, \ \ \  \tau=t.$$
 Then the free boundaries $x=a(t)$ and $x=b(t)$ become $\xi=0$ and $\xi=\int_{a(t)}^{b(t)}m(y,t)dy=
 \int_{a}^{b}m_{0}(y)dy$ by the conservation of mass, where $\int_{a}^{b}m_{0}(y)dy$ is the total  liquid mass
 initially. We normalize $\int_{a}^{b}m_{0}(y)dy$ to $1$.

 Hence in the Lagrangian coordinates, the two free boundary problems (1.6), (1.10),
(1.11) and (1.6), (1.12), (1.13) become
$$
\left \{
\begin{array}{l}
n_{\tau}+  (nm)u_{\xi} = 0,\\[3mm]
m_{\tau}+  m^{2}u_{\xi} = 0, \\[3mm]
u_{\tau}+(P(n,m))_{\xi}=(\varepsilon(n,m)m u_{\xi})_{\xi},
\end{array}
\right. \eqno(2.1)
$$
with the boundary conditions (corresponding to the initial  masses
connect to vacuum discontinuously)
     $$
     P(n,m)=E(n,m)u_{\xi}, \ \ \ {\rm at} \ \ \  \xi=0,1, \ \ \tau \geq 0, \eqno(2.2)
     $$
and the initial data
$$n(\xi,0)=n_{0}(\xi)>0,\ \ m(\xi,0)=m_{0}(\xi)>0,\ \ u(\xi,0)=u_{0}(\xi),
\ \ \ \xi \in [0,1], \eqno(2.3)$$

\noindent or with the boundary conditions (corresponding to the
initial masses connect to vacuum continuously)
     $$n(0,\tau)=n(1,\tau)=0, \ \ \ m(0,\tau)=m(1,\tau)=0, \ \ \ \tau \geq 0,\eqno(2.4)$$
and the initial data
$$n(\xi,0)=n_{0}(\xi)>0,\ \ m(\xi,0)=m_{0}(\xi)>0,\ \ u(\xi,0)=u_{0}(\xi),
\ \ \ \xi \in (0,1),\eqno(2.5)$$ and
$n_{0}(0)=n_{0}(1)=m_{0}(0)=m_{0}(1)=0$.

Here
$$P(n,m)= \left (\frac{n}{\rho_{l}-m}\right)^{\gamma}, \ \ \
 \varepsilon(n,m)m=\frac{n^{\beta}m}{(\rho_{l}-m)^{\beta+1}}= E(n,m), \ \ \ \beta>0.\eqno(2.6)$$
Here we have assumed $A=B=1$ in (2.6) for simplicity.

In the following, we replace the coordinates $(\xi, \tau)$ by
$(x,t)$. Introduce the variables (cf.
\cite{Evje-Karlsen,Evje-Flatten-Friis}):
$$c=\frac{n}{m},\ \ \ \ \ \
Q(m)=\frac{m}{\rho_{l}-m}=\frac{\alpha_{l}}{1-\alpha_{l}}\geq 0.$$
Form the first two equations of (2.1), we get
$$c_{t}=\frac{n_{t}}{m}-\frac{n}{m^{2}}m_{t}=-\frac{mnu_{x}}{m}+\frac{nm^{2}}{m^{2}}u_{x}=0,$$
and
$$\arraycolsep=1.5pt
\begin{array}{rcl}
Q(m)_{t}=\D\left(\frac{m}{\rho_{l}-m}\right)_{t}&=&
\D \left(\frac{1}{\rho_{l}-m}+\frac{m}{(\rho_{l}-m)^{2}}\right)m_{t}\\[3mm]
&=&\D \frac{\rho_{l}}{(\rho_{l}-m)^{2}}m_{t}=-\D\frac{\rho_{l}m^{2}}{(\rho_{l}-m)^{2}}u_{x}\\[3mm]
&=&-\rho_{l}Q(m)^{2}u_{x}.
\end{array}
$$
Then we can rewrite the initial boundary problems (2.1), (2.2),
(2.3) and (2.1), (2.4), (2.5) into the following forms:
$$
\left \{
\begin{array}{l}
\partial_{t}c=0,\\[3mm]
\partial_{t}Q(m)+\rho_{l}Q(m)^{2}\partial_{x}u = 0,  \\[3mm]
\partial_{t}u+\partial_{x}(P(c,m))=\partial_{x}(E(c,m)\partial_{x}u),
\end{array}
\right. \eqno(2.7)
$$
with the boundary conditions (corresponding to the initial  masses
connect to vacuum discontinuously)
     $$
     P(c,m)=E(c,m)u_{x},\ \ \ \ {\rm at} \ \ \  x=0,1, \ \ t\geq 0,\eqno(2.8)
     $$
and the initial data
$$c(x,0)=c_{0}(x)>0,\ Q(m)(x,0)=Q(m_{0})(x)>0,\ \ u(x,0)=u_{0}(x),
\ \ \ x \in [0,1],\eqno(2.9)$$

\noindent or with the boundary conditions (corresponding to the
initial masses connect to vacuum continuously)
     $$c(0,t)=c(1,t)=0, \ \ \ Q(m)(0,t)=Q(m)(1,t)=0, \ \ \ t\geq 0, \eqno(2.10)$$
and the initial data
$$c(x,0)=c_{0}(x)> 0,\ Q(m)(x,0)=Q(m_{0})(x)> 0,\ u(x,0)=u_{0}(x),
\ \ \ x \in (0,1), \eqno(2.11)$$ and
$c_{0}(0)=c_{0}(1)=Q(m_{0})(0)=Q(m_{0})(1)=0$.

Here
$$P(c,m)= \left (\frac{n}{\rho_{l}-m}\right)^{\gamma}=c^{\gamma}Q(m)^{\gamma},\ \ \
 E(c,m)= m\varepsilon(n,m)=c^{\beta}Q(m)^{\beta+1}, \ \ \ \beta>0.\eqno(2.12)$$

\vspace{4mm} \noindent Throughout this paper, our assumptions on the
initial data and $\beta,\ \gamma$ are stated as follows:

\vspace{2mm}

$(A_1)$ \ \ $\inf\limits_{x\in [0,1]}n_0(x)>0, \ \ \
\sup\limits_{x\in [0,1]}n_0(x)<\infty, \ \ \ \inf\limits_{x\in
[0,1]}m_0(x)>0, \ \ \ \sup\limits_{x\in [0,1]}m_0(x)<\rho_l;$

\vspace{2mm}

$(A_{1})'$ \ \ There are positive constants $K_{1},\ K_{2},\ K_{3} $
and $K_{4}$ such that $K_{1}\phi(x)^{\frac{\alpha}{2}} \leq
m_{0}(x)\leq K_{2}\phi(x)^{\frac{\alpha}{2}}<\rho_{l}$,
$K_{3}\phi(x)^{\alpha} \leq  n_{0}(x)\leq K_{4}\phi(x)^{\alpha}$,
where $\phi(x)=x(1-x)$, $0<\alpha<1$. In particular, this implies
that there exist positive constants $C_{1}, \ C_{2}$, shuch that
$C_{1}\phi(x)^{\frac{\alpha}{2}} \leq
c_{0}(x)=\frac{n_{0}(x)}{m_{0}(x)}\leq
C_{2}\phi(x)^{\frac{\alpha}{2}}$;

\vspace{2mm}

$(A_{2})\ \ u_{0}(x) \in L^{2n}([0,1])$ for any given positive
integer $n$ satisfying $n \geq  \frac{2\gamma+\beta}{2\beta}$;

\vspace{2mm}

$(A_{3})\ \ ((c_{0}Q(m_{0}))^{\beta})_{x} \in L^{2}([0,1])$;

\vspace{2mm}

$(A_4)$ \ \ $\beta >0,\gamma \geq 1+\beta $.

\vspace{4mm}

Now we give the following definition of weak solution:

\vspace{2mm}

\noindent{\bf Definition 2.1 (Weak solution).} \ \ We call $(n(x,t),
m(x,t), u(x,t)) $  a global weak solution to the initial boundary
value problems (2.1), (2.2), (2.3) or (2.1), (2.4) (2.5), if the
following estimates hold for any $t>0$,
$$n,\ m,\ u \in L^{\infty}([0,1]\times[0,+\infty) )\cap
C^{1}([0,+\infty);H^{1}([0,1])),$$
$$E(n,m)u_{x} \in L^{\infty}([0,1]\times[0,+\infty) )\cap
C^{\frac{1}{2}}([0,+\infty);L^{2}([0,1])),$$
$$
0 \leq n(x,t)< \rho_{l}\sup_{x\in [0,1]}c_{0},
$$
and
$$
0 \leq m(x,t)< \rho_{l}.
$$
Furthermore, the following equations hold:
$$
n_{t}+mnu_{x}=0,\ \ \ \ \  m_{t}+m^{2}u_{x}=0,
$$
$$
(n,m)(x,0)=(n_{0}(x),m_{0}(x)), \ \ \ \textrm{for} \ \  {\rm a.e.} \
\  x \in [0,1] \ \ \textrm{and} \ \ \textrm{any} \ \  t \geq 0,
$$
and
$$
\int_0^{\infty}\int_{0}^{1}(
u\varphi_{t}+(P(n,m)-E(n,m)u_x)\varphi_x)dxdt
+\int_{0}^{1}u_0(x)\varphi(x,0)dx=0,
$$
for any test functions $\varphi \in C_0^{\infty}(\Omega)$ with
$\Omega=\{(x,t): \ 0 \leq x\leq 1, \ t\geq 0\}$.

\vspace{4mm}

In what follows, we always use $C$ (and $C_n$) to denote a generic
positive constant depending only on the initial data (and the given
positive integer $n$), but independent of $t$.

\vspace{4mm}

 We now state the main theorems in this paper as
follows:
 \vspace{2mm}

 \noindent{\bf Theorem 2.2 (The asymptotic behavior of the mass functions).}
 \ \ Under the assumptions $(A_{1})$, $(A_{2})$, $(A_{3})$ and
 $(A_{4})$ (corresponding to the boundary conditions $(2.8)$) or
 $(A_{1})'$, $(A_{2})$, $(A_{3})$ and
 $(A_{4})$ (corresponding to the boundary conditions $(2.10)$),
  let $(n(x, t), \ m(x,t),$ $u(x,t))$ be a global weak solution to the initial boundary
value problem (2.1), (2.2), (2.3) or (2.1), (2.4) (2.5). Then we
have the following asymptotic behavior of the mass functions $n(x,
t),\ m(x,t)$
$$
\lim\limits_{t \to \infty}\sup\limits_{x\in [0, 1]}n(x, t)=0,
\eqno(2.13)
$$
$$
\lim\limits_{t \to \infty}\sup\limits_{x\in [0, 1]}m(x,
t)=0. \eqno(2.14)
$$

\vspace{4mm}
 Furthermore, we can get the decay rate estimates of the
mass functions $n(x,t)$,\ $ m(x,t)$ as follows:

\vspace{2mm}

\noindent{\bf Theorem 2.3 (The decay rate of the mass functions).} \
\ Under the assumptions of Theorem 2.2, let $(n(x, t), \ m(x,t),$
$u(x,t))$ be a global weak solution to the initial boundary value
problem (2.1), (2.2), (2.3) or (2.1), (2.4) (2.5). Then the
following decay rate estimates on the mass functions $n(x, t),\
m(x,t)$ hold:

 (i) Under the
boundary conditions $(2.8)$, if $ 0< \beta <1$ or $\beta>1,\
\frac{\gamma-1}{\gamma-\beta}>2$, then
 $$
 n(x,t),\ m(x,t)\leq C(1+t)^{-\frac{\theta}{\gamma-1+2\beta}}, \eqno(2.15)
$$
for  any $x \in [0,1]$.

If $\beta=1$ or $\beta>1,\ \frac{\gamma-1}{\gamma-\beta}\leq2$, we
have
$$
 n(x,t),\ m(x,t)\leq C(1+t)^{-\frac{\theta}{\gamma-1+2\beta}}(\ln(1+t))^{\frac{1}{\gamma-1+2\beta}}, \eqno(2.16)
$$
for any  $ x \in [0,1]$.

(ii) Under the boundary conditions $(2.10)$, if $ 0< \beta <1$ or
$\beta>1,\ \frac{\gamma-1}{\gamma-\beta} >2$, then
$$
 n(x,t),\ m(x,t) \leq C(1+t)^{-\frac{\theta}{\gamma-1+4\beta}}, \eqno(2.17)
$$
for any  $ x \in [0,1]$.

If $\beta=1$ or $\beta>1,\ \frac{\gamma-1}{\gamma-\beta} \leq2$, we
have
 $$n(x,t),\ m(x,t)\leq C(1+t)^{-\frac{\theta}{\gamma-1+4\beta}}(\ln(1+t))^{\frac{1}{\gamma-1+4\beta}}, \eqno(2.18)
$$
for  any $ x \in [0,1]$, where $\theta$ is defined by Lemma 4.1.

\section{A \emph{priori} estimates and the asymptotic behavior of the mass functions }

In this section, we will give some useful uniform {\it a priori}
estimates of the solutions with respect to the time $t$. Then we
study the asymptotic behavior of the mass functions $n(x,t)$ and
$m(x,t)$ by using these uniform {\it  a priori} estimates.

\subsection{Uniform \emph{a priori} estimates}

\noindent{\bf Lemma 3.1  (Some identities)}.

$$
c(x,t)=c_{0}(x),\eqno(3.1)
$$
$$
\displaystyle \frac{d}{dt} \int_{0}^{x} u(y,t)dy = -\frac{d}{dt}
\int_{x}^{1}u(y,t)dy,\eqno(3.2)$$
$$
\displaystyle
(c^{\beta}Q(m)^{\beta+1}u_{x})(x,t) =
(c^{\gamma}Q(m)^{\gamma})(x,t)+\int_{0}^{x}u_{t}(y,t)dy
 =(c^{\gamma}Q(m)^{\gamma})(x,t)-\int_{x}^{1}u_{t}(y,t)dy,\eqno(3.3)
 $$
$$\displaystyle \frac{1}{\beta
\rho_{l}}(c^{\beta}Q(m)^{\beta})(x,t)+\int_{0}^{t}c^{\gamma}Q(m)^{\gamma}(x,s)ds=
\frac{1}{\beta\rho_{l}} c
_{0}^{\beta}Q(m_{0})^{\beta}-\int_{0}^{x}\int_{0}^{t}u_{t}(y,s)dyds.\eqno(3.4)
 $$

 \vspace{2mm}

\noindent {\it Proof.} \ \ These identities can be obtained directly
from $(2.7)$.

  \vspace{4mm}

\noindent{\bf Lemma 3.2 (Basic energy estimate).} \ \ Under the
conditions in Theorem 2.2, the following ener\/gy estimate holds:
$$
\arraycolsep=1.5pt
\begin{array}[b]{rl} &\displaystyle \int_{0}^{1}\left( \frac{1}{2}
u^{2}+\frac{c^{\gamma}}{\rho_{l}(\gamma-1)}Q(m)^{\gamma-1}(x,t)\right)dx
+\int_{0}^{t}\int^{1}_{0}c^{\beta}Q(m)^{\beta+1}u_{x}^{2}dx ds \\[3mm]
=& \displaystyle
 \int_{0}^{1}\left( \frac{1}{2} u_{0}^{2} + \frac{c_{0}^{\gamma}}{\rho_{l}(\gamma-1)} Q(m_{0})^{\gamma-1}\right)dx \\[3mm]
\leq & \displaystyle C.
\end{array}\eqno(3.5)
$$

\noindent {\it Proof.}\ \ Multiplying the second and the third
equations of (2.7) by $c^{\gamma}Q(m)^{\gamma-2}$ and $u$, and
integrating the resulting equations with respect to $x$ over [0,1],
we get
$$
\begin{array}[b]{rl}
\arraycolsep=1.5pt
& \displaystyle \frac{d}{dt}\int^{1}_{0}\left(\frac{u^{2}}{2}
+\frac{c^{\gamma}}{\rho_{l}(\gamma-1)}Q(m)^{\gamma-1}\right)dx+(u P(c,m))\Big|_{0}^{1}\\[3mm]
=&(E(c,m)u_{x}u) \Big|_{0}^{1}-
\displaystyle\int_{0}^{1}E(c,m)u_{x}^{2}dx.
\end{array}\eqno(3.6)
$$
Using the boundary conditions $(2.8)$ or $(2.10)$ we get
$$
\displaystyle
\frac{d}{dt}\int^{1}_{0}\left(\frac{u^{2}}{2}+\frac{c^{\gamma}}{\rho_{l}(\gamma-1)}Q(m)^{\gamma-1}\right)dx
+\int_{0}^{1}c^{\beta}Q(m)^{\beta+1}u_{x}^{2}dx=0.  \eqno(3.7)
$$
Then integrating it with respect to $t$ over $[0,t]$, we get
$(3.5)$.

The proof of Lemma 3.2 is completed.

\vspace{4mm}

\noindent { \bf Lemma 3.3 (The Uniform upper bound for the
$cQ(m)$).}  Under the conditions of Theorem $2.2$, we have  for any
$x \in[0,1],\ t>0$,
$$ 0   \leq (cQ(m))(x,t)\leq C .\eqno(3.8)$$

\vspace{3mm}

\noindent {\it Proof.}\ \ By (3.4), we have
$$
\arraycolsep=1.5pt
\begin{array}{rl}
\D \frac{1}{\beta
\rho_{l}}(c^{\beta}Q(m)^{\beta})(x,t)+\int_{0}^{t}c^{\gamma}Q(m)^{\gamma}(x,s)ds
=&\D \frac{1}{\beta\rho_{l}}
c_{0}^{\beta}Q(m_{0})^{\beta}-\int_{0}^{x}(u(y,t)-u_{0}(y))dy\\[3mm]
\leq & \D C+\int_{0}^{1}|u(y,t)|dy +\int_{0}^{1}|u_{0}(y)|dy\\[3mm]
\leq & \D C+ \left(\int_{0}^{1}u^{2}dy
\right)^{\frac{1}{2}}+\left(\int_{0}^{1}u_{0}^{2}dy\right)^{\frac{1}{2}}\\[3mm]
\leq & C.
\end{array}
$$
Here we have used Lemma $3.2$ and the assumptions $(A_1), \ (A_1)',
\ (A_2)$.

The proof of Lemma 3.3 is completed.

\vspace{4mm}

\noindent{\textbf{Corollary 3.4.}} \ \ For any $x\in[0,1]$ and
$t>0$,
$$
\int_{0}^{t}c^{\gamma}Q(m)^{\gamma}(x,s)ds \leq C. \eqno(3.9)
$$

\noindent{\bf Lemma 3.5.} \ \ Under the boundary conditions (2.8),
we have for any $t>0$,
$$
Q(m)(d,t)=Q(m_{0})(d)\left(\frac{1}{(\gamma-\beta)\rho_{l}c_{0}^{\gamma-\beta}(d)Q(m_{0})^{\gamma-\beta}t+1}
\right)^{\frac{1}{\gamma-\beta}},\ \ \ \ d=0, \ 1.
\eqno(3.10)
$$

\vspace{3mm}

\noindent {\it Proof.}\ \ By $(3.2)$, we have
$$
\frac{d}{dt} \int_{0}^{1} u(y,t)dy = 0.
$$
By taking $x=1$ or $x=0$ in $(3.4)$, we have
$$
\frac{1}{\beta
\rho_{l}}c_{0}^{\beta}(d)Q(m)^{\beta}(d,t)+\int_{0}^{t}c_{0}^{\gamma}(d)Q(m)^{\gamma}(d,s)ds=
\frac{1}{\beta\rho_{l}} c _{0}^{\beta}(d)Q(m_{0})^{\beta}(d),\ \ \
d=0,1. \eqno(3.11)
$$
Since $\gamma > \beta$, the integral equation $(3.11)$ yields
$$
\frac{1}{\rho_{l}}c_{0}^{\beta}(d)Q(m)^{\beta-1}(d,t)Q(m)(d,t)_{t}+c_{0}^{\gamma}(d)Q(m)^{\gamma}(d,t)=0,
\eqno(3.12)
$$
which implies (3.10)  by solving the ordinary  differential equation
(3.12).

The proof of Lemma $3.5$ is completed.

\vspace{4mm}

\noindent{\textbf{Corollary 3.6.}} \ \ There exist positive
constants $C_{1}$ and $C_{2}$ such that for any $t>0$,
$$
C_{1}(1+t)^{-\frac{1}{\gamma-\beta}} \leq Q(m)(d,t)\leq C_{2}(1+t)^{-\frac{1}{\gamma-\beta}}.
$$

\vspace{3mm}

 \noindent{\bf Lemma 3.7.} \ \ For any  positive integer
$n$ in $(A_2)$, we have for any $t>0$
$$
\int_{0}^{1}u^{2n}dx+n(2n-1)\int_{0}^{t}\int_{0}^{1}c^{\beta}Q(m)^{\beta+1}u^{2n-2}u_{x}^{2}dxds
\leq C_{n}, \eqno(3.13)
$$
where $C_{n}$ is a positive constant
depending  on $n$, but independent of $t$.

\vspace{3mm}

\noindent {\it Proof.}\ \ Multiplying the third equation of $(2.7)$
by $2nu^{2n-1}$ and integrating the resulting equation with respect
to $x$ over $[0,1]$, we have
$$
\arraycolsep=1.5pt
\begin{array}[b]{rl}
& \displaystyle \frac{d}{dt}\int^{1}_{0}u^{2n}dx+2n (u^{2n-1}
P(c,m))\Big|_{0}^{1}
-2n(2n-1)\int_{0}^{1} c^{\gamma}Q(m)^{\gamma}u^{2n-2}u_{x}dx\\[3mm]
=&2n(u^{2n-1}E(c,m)u_{x})\Big|_{0}^{1}
-2n(2n-1)\displaystyle\int_{0}^{1}c^{\beta}Q(m)^{\beta+1}u^{2n-2}u_{x}^{2}dx.
\end{array}
\eqno(3.14)
$$
Using the boundary conditions $(2.8)$ or $(2.10)$, we have
$$
\displaystyle \frac{d}{dt}\int^{1}_{0}u^{2n}dx+
2n(2n-1)\displaystyle\int_{0}^{1}c^{\beta}Q(m)^{\beta+1}u^{2n-2}u_{x}^{2}dx=2n(2n-1)\int_{0}^{1}
c^{\gamma}Q(m)^{\gamma}u^{2n-2}u_{x}dx . \eqno(3.15)
$$ Integrating
$(3.15)$ with respect to $t$ over $[0,t]$, we get
$$
\begin{array}[b]{rl}
&\D \int^{1}_{0}u^{2n}dx+
2n(2n-1)\displaystyle\int_{0}^{t}\int_{0}^{1}c^{\beta}Q(m)^{\beta+1}u^{2n-2}u_{x}^{2}dxds\\[3mm]
=& \D \int^{1}_{0}u_{0}^{2n}dx+2n(2n-1)\int_{0}^{t}\int_{0}^{1}
c^{\gamma}Q(m)^{\gamma}u^{2n-2}u_{x}dxds.
\end{array}
\eqno(3.16)
$$
Applying Cauchy-Schwarz inequality to the last term
in $(3.16)$ yields
$$
\arraycolsep=1.5pt
\begin{array}{rl}
&\D \int^{1}_{0}u^{2n}dx+
n(2n-1)\displaystyle\int_{0}^{t}\int_{0}^{1}c^{\beta}Q(m)^{\beta+1}u^{2n-2}u_{x}^{2}dxds\\[3mm]
\leq & \D \int^{1}_{0}u_{0}^{2n}dx+n(2n-1)\int_{0}^{t}\int_{0}^{1}
c^{2\gamma-\beta}Q(m)^{2\gamma-\beta-1}u^{2n-2}dxds.
\end{array}
\eqno(3.17)
$$
Now we estimate the last term on the right-hand side in $(3.17)$ as
follows:
$$
\arraycolsep=1.5pt
\begin{array}{rl}
&n(2n-1)\D\int_{0}^{t} \int_{0}^{1}c^{2\gamma-\beta}Q(m)^{2\gamma-\beta-1}u^{2n-2}dxds\\[3mm]
=& n(2n-1)\D \int_{0}^{t}\int_{0}^{1}c^{ \frac{\gamma}{n}+\gamma
-\beta}Q(m)^{\frac{\gamma}{n}
+\gamma-\beta-1}c^{\frac{n-1}{n}\gamma}Q(m)^{\frac{n-1}{n}\gamma}u^{2n-2}dxds\\[3mm]
\leq &(2n-1)\D \int_{0}^{t}\int_{0}^{1}c^{ \gamma+n(\gamma
-\beta)}Q(m)^{\gamma+n(\gamma-\beta-1)}dxds+(n-1)(2n-1)\int_{0}^{t}\int_{0}^{1}c^{\gamma}Q(m)^{\gamma}u^{2n}dxds\\[3mm]
=&(2n-1)\D \int_{0}^{t}\int_{0}^{1}(cQ(m))^{n(\gamma
-\beta-1)}c^{n}(cQ(m))^{\gamma}dxds+(n-1)(2n-1)\int_{0}^{t}\int_{0}^{1}c^{\gamma}Q(m)^{\gamma}u^{2n}dxds\\[3mm]
\leq & C \D \int_{0}^{t}\max \limits _{[0,1]}(cQ(m))^{\gamma}ds
+C\int_{0}^{t}\max
\limits_{[0,1]}(cQ(m))^{\gamma}\int_{0}^{1}u^{2n}dxds\\
\leq &C+\D
C\int_{0}^{t}\max\limits_{[0,1]}(cQ(m))^{\gamma}\int_{0}^{1}u^{2n}dxds.
\end{array}
\eqno(3.18)
$$
Here we  have used the Young inequality $ ab \leq
\frac{a^{p}}{p}+\frac{b^{q}}{q}$, where $\frac{1}{p}+\frac{1}{q}=1,\
p,q>1,\ a,\ b \geq 0$, the assumption $(A_{4})$, $(3.8)$ and
$(3.9)$.

Substituting (3.18) into (3.17), we have
$$
\begin{array}[b]{rl}
&\D \int^{1}_{0}u^{2n}dx+
n(2n-1)\displaystyle\int_{0}^{t}\int_{0}^{1}c^{\beta}Q(m)^{\beta+1}u^{2n-2}u_{x}^{2}dxds\\[3mm]
\leq & C+\D
C\int_{0}^{t}\max\limits_{[0,1]}(cQ(m))^{\gamma}\int_{0}^{1}u^{2n}dxds.
\end{array}
\eqno(3.19)
$$
By (3.19), (3.9) and Gronwall's inequality, we have
$$
\int_0^1u^{2n}dx\leq C_n, \eqno(3.20)
$$
where $C_n$ is a positive constant depending on $n$, but independent
of $t$.

(3.19) and (3.20) show that (3.13) holds and this completes the
proof of Lemma 3.7.

\vspace{4mm}

\noindent {\bf Lemma 3.8.} \ \ We have the following uniform
estimate on the derivative of the function $cQ(m)$,
$$
\D \int_{0}^{1} \left((cQ(m))^{\beta} \right)_{x}^{2}dx+
\int_{0}^{1}\int_{0}^{t}\left((cQ(m))^{\frac{\beta+\gamma}{2}}
\right)_{x}^{2}dsdx \leq C. \eqno(3.21)
$$

\vspace{3mm}

\noindent {\it Proof.}\ \ From $(2.7)$, we have
$$
\arraycolsep=1.5pt
\begin{array}[b]{rcl}
\D ((cQ(m))^{\beta})_{xt}&=& \D(\beta(cQ(m))^{\beta-1}(cQ(m))_{t})_{x}\\[3mm]
                      &=& \D(\beta c^{\beta}Q(m)^{\beta-1}Q(m)_{t})_{x}\\[3mm]
                      &=&\D-\beta \rho_{l}(c^{\beta}Q(m)^{\beta+1}u_{x})_{x}\\[3mm]
&=&\D-\beta \rho_{l}(u_{t}+P(c,m)_{x}).
\end{array}
\eqno(3.22)
$$
Multiplying (3.22) by $((cQ(m))^{\beta})_{x}$ ,
integrating the resulting equation  over $[0,1]\times[0,t]$ and
integrating by parts, we get
$$
\arraycolsep=1.5pt
\begin{array}{rcl}
\D \frac{1}{2} \int_{0}^{1}\left((cQ(m))^{\beta} \right)_{x}^{2}dx&=
& \D \frac{1}{2} \int_{0}^{1}\left((c_{0}Q(m_{0}))^{\beta} \right)_{x}^{2}dx
-\beta \rho_{l}\int_{0}^{1}u\left((cQ(m))^{\beta} \right)_{x}dx\\[3mm]
&&+\D \beta \rho_{l}\int_{0}^{1}u_{0}\left((c_{0}Q(m_{0}))^{\beta} \right)_{x}dx
+\D \beta \rho_{l}\int_{0}^{1}\int_{0}^{t}u\left((cQ(m))^{\beta} \right)_{xt}dsdx\\[3mm]
&&-\D\frac{4\gamma \beta^{2}\rho_{l}}{(\beta+\gamma)^{2}}
\int_{0}^{1}\int_{0}^{t}\left((cQ(m))^{\frac{\beta+\gamma}{2}}\right)_{x}^{2}dsdx.
\end{array}
\eqno(3.23)$$ Substituting (3.22) into (3.23), we have
$$
\arraycolsep=1.5pt
\begin{array}[b]{rl}
& \D \frac{1}{2} \int_{0}^{1}\left((cQ(m))^{\beta} \right)_{x}^{2}dx+
\D\frac{4\gamma \beta^{2}\rho_{l}}{(\beta+\gamma)^{2}}
\int_{0}^{1}\int_{0}^{t}\left((cQ(m))^{\frac{\beta+\gamma}{2}}\right)_{x}^{2}dsdx\\[3mm]
=&\D \frac{1}{2} \int_{0}^{1}\left((c_{0}Q(m_{0}))^{\beta} \right)_{x}^{2}dx
-\beta \rho_{l}\int_{0}^{1}u\left((cQ(m))^{\beta} \right)_{x}dx\\[3mm]
&+\D \beta \rho_{l}\int_{0}^{1}u_{0}\left((c_{0}Q(m_{0}))^{\beta} \right)_{x}dx
-(\beta\rho_{l})^{2}\int_{0}^{1}\int_{0}^{t}uu_{t}dsdx\\[3mm]
&-\D(\beta\rho_{l})^{2}\int_{0}^{1}\int_{0}^{t}u((cQ(m))^{\gamma})_{x}dsdx\\[3mm]
=&\D \frac{1}{2} \int_{0}^{1}\left((c_{0}Q(m_{0}))^{\beta} \right)_{x}^{2}dx
-\beta \rho_{l}\int_{0}^{1}u\left((cQ(m))^{\beta} \right)_{x}dx\\[3mm]
&+\D \beta \rho_{l}\int_{0}^{1}u_{0}\left((c_{0}Q(m_{0}))^{\beta} \right)_{x}dx
-\frac{(\beta \rho_{l})^{2} }{2}\int_{0}^{1}u^{2}dx+\frac{(\beta \rho_{l})^{2} }{2}\int_{0}^{1}u_{0}^{2}dx\\[3mm]
&\D -(\beta \rho_{l})^{2}\int_{0}^{t}\{(cQ(m))^{\gamma}(1,s)u(1,s)-(cQ(m))^{\gamma}(0,s)u(0,s)\}ds\\[3mm]
&\D+(\beta\rho_{l})^{2}\int_{0}^{1}\int_{0}^{t}(cQ(m))^{\gamma}u_{x}dsdx
=\sum \limits_{i=1}^{i=7}J_{i}.
\end{array}
\eqno(3.24)
$$
Now we estimate $J_1$--$J_7$ as follows:

First, by the assumptions $(A_{1})$-$(A_{3})$, or
$(A_{1})'$-$(A_{3})$, Lemma 3.2, Lemma 3.3, Corollary 3.4 and
Cauchy-Schwarz inequality, we have
$$
\left\{
\begin{array}{l}
 \displaystyle J_1\leq C, \\ [3mm]
 \displaystyle J_2\leq \frac{1}{4}\int_{0}^{1}\left((cQ(m))^{\beta} \right)_{x}^{2}dx
 +C\int_{0}^{1}u^{2}dx\leq C+\frac{1}{4}\int_{0}^{1}\left((cQ(m))^{\beta} \right)_{x}^{2}dx,\\ [3mm]
\displaystyle J_3 \leq  C\int_{0}^{1}\left((c_{0}Q(m_{0}))^{\beta}
\right)_{x}^{2}dx +C\int_{0}^{1}u_{0}^{2}dx \leq C,\\ [3mm]
 \displaystyle J_4\leq C, \\ [3mm]
\displaystyle J_5\leq C, \\ [3mm]
\arraycolsep=1.5pt
\begin{array}{rl}
\displaystyle J_7\leq  & \displaystyle
C\int_{0}^{t}\int_{0}^{1}c_0(x)c^{\beta}Q(m)^{\beta+1}u_{x}^{2}dxds
+C\int_{0}^{t}\int_{0}^{1}(cQ(m))^{2\gamma-\beta-1}dxds\\ [3mm] \leq
& \displaystyle C\max c_0(x)
\int_{0}^{t}\int_{0}^{1}c^{\beta}Q(m)^{\beta+1}u_{x}^{2}dxds
+C\max(cQ(m))^{\gamma-\beta-1}\int_{0}^{t}\max(cQ(m))^{\gamma}ds \\[3mm]
\leq & C.
\end{array}
 \end{array}
\right. \eqno(3.25)
$$
Now we estimate $J_{6}$, which is divided into two cases:

Case 1. When the initial masses connect to vacuum continuously
(corresponding to the boundary condition $(2.10)$), we have
$J_{6}=0$.

\vspace{2mm}

Case 2. When the initial masses connect to vacuum discontinuously
(corresponding to the boundary condition $(2.8)$), we have by Young
 inequality and Lemma 3.5
$$
\arraycolsep=1.5pt
\begin{array}[b]{rl}
\displaystyle J_6=  &
 \displaystyle-(\beta\rho_{l})^2\int_0^{t}(cQ(m))^{\gamma-\beta}(1,
s)\left((cQ(m))^{\beta}(1, s)u(1, s)\right)ds \\[3mm]
&\displaystyle+(\beta\rho_{l})^2\int_0^{t}(cQ(m))^{\gamma-\beta}(0,s)\left(cQ(m))^{\beta}(0,
s)u(0, s)\right)ds \\ [3mm]
 \leq & \displaystyle
C\int_0^{t}\left\{\left|(cQ(m))^{n\beta}(1, s)u^n(1,s)\right|
+\left|(cQ(m))^{n\beta}(0, s)u^n(0, s)\right|\right\}ds \\[3mm]
& \displaystyle
+C\int_0^{t}\left\{(cQ(m))^{(\gamma-\beta)\frac{n}{n-1}}(1, s)
+(cQ(m))^{(\gamma-\beta)\frac{n}{n-1}}(0, s)\right\}ds \\ [3mm]
 \leq
& \displaystyle C+C\int_0^{t}\left|\left|(cQ(m))^{n\beta}(\cdot,
s)u^n(\cdot, s)\right|\right|_{L^{\infty}([0, 1])}ds.
\end{array}\eqno(3.26)
$$
Substituting (3.25) and (3.26) into (3.24), we have for Case 1 and
Case 2
$$
\D\frac{1}{2} \int_{0}^{1} \left((cQ(m))^{\beta} \right)_{x}^{2}dx+
\frac{4\gamma\beta^{2}\rho_{l}}{(\beta+\gamma)^{2}}
\int_{0}^{1}\int_{0}^{t}\left((cQ(m))^{\frac{\beta+\gamma}{2}}
\right)_{x}^{2}dsdx \leq C +J_8,\eqno(3.27)
$$
where
$$
J_8=\displaystyle C\int_0^{t}\left|\left|(cQ(m))^{n\beta}(\cdot,
s)u^n(\cdot, s)\right|\right|_{L^{\infty}([0, 1])}ds.
$$
By the embedding theorem $W^{1,1}([0, 1])\hookrightarrow
L^{\infty}([0, 1])$, we have
$$
\arraycolsep=1.5pt
\begin{array}[b]{rcl}
\displaystyle J_8&\leq&\D
C\int_{0}^{t}\int_{0}^{1}(cQ(m))^{n\beta}u^{n}(x,s)dxds+
\D C \int_{0}^{t}\int_{0}^{1}\left((cQ(m))^{n\beta}u^{n}(x,s)\right)_{x}dxds\\[3mm]
&\leq &C \D \int_{0}^{t}\int_{0}^{1}(cQ(m))^{\gamma}u^{2n}(x,s)dxds
+C\D \int_{0}^{t}\int_{0}^{1}(cQ(m))^{2n\beta-\gamma}dxds\\[3mm]
&&\D+ \int_{0}^{t}\int_{0}^{1}n\beta(cQ(m))^{n\beta-1}(cQ(m))_{x}u^{n}(x,s)dxds
\D+ \int_{0}^{t}\int_{0}^{1}n(cQ(m))^{n\beta}u^{n-1}u_{x}dxds\\[3mm]
&\leq & C\D \int_{0}^{t}\max\limits_{[0,1]}(cQ(m))^{\gamma}\left(\int_{0}^{1}u^{2n}dx\right)ds
+\D C\int_{0}^{t}\max\limits_{[0,1]}(cQ(m))^{\gamma}\int_{0}^{1}(cQ(m))^{2n\beta-2\gamma}dxds\\[3mm]
&&\D+C
\int_{0}^{t}\int_{0}^{1}(cQ(m))^{2n\beta-\gamma-\beta}u^{2n}dxds
+\frac{2\gamma\beta^{2}\rho_{l}}{(\beta+\gamma)^{2}}\int_{0}^{1}\int_{0}^{t}\left((cQ(m))^{\frac{\beta+\gamma}{2}}
\right)_{x}^{2}dsdx\\[3mm]
&&\D+C \int_{0}^{t}\int_{0}^{1}(cQ(m))^{\beta+1}u^{2n-2}u_{x}^{2}dxds+C
\int_{0}^{t}\int_{0}^{1}(cQ(m))^{2n\beta-\beta-1}dxds\\[3mm]
& \leq & C+ \D \max(cQ(m))^{2n\beta-2\gamma-\beta}
\int_{0}^{t}\max\limits_{[0,1]}(cQ(m))^{\gamma}\left(\int_{0}^{1}u^{2n}dx\right)ds\\[3mm]
&&\D+\frac{2\gamma\beta^{2}\rho_{l}}{(\beta+\gamma)^{2}}\int_{0}^{1}\int_{0}^{t}\left((cQ(m))^{\frac{\beta+\gamma}{2}}
\right)_{x}^{2}dsdx+C\D \max(c_0(x))\int_{0}^{t}\int_{0}^{1}c^{\beta}Q(m)^{\beta+1}u^{2n-2}u_{x}^{2}dxds\\[3mm]
&& +C\D
\max(cQ(m))^{2n\beta-\beta-1-\gamma}\int_{0}^{t}\max\limits_{[0,1]}(cQ(m))^{\gamma}ds\\[3mm]
&\leq & C+\D\frac{2\gamma\beta^{2}\rho_{l}}{(\beta+\gamma)^{2}}
\int_{0}^{1}\int_{0}^{t}\left((cQ(m))^{\frac{\beta+\gamma}{2}}
\right)_{x}^{2}dsdx.\\[3mm]
\end{array}
\eqno(3.28)
$$
Here we have used Lemma 3.3, Corollary 3.4, Lemma 3.7, and $n\geq
\frac{2\gamma+\beta}{2\beta}$.

Substituting (3.28) into (3.27), we get (3.21). This proves Lemma
3.8. \vspace{2mm}

\subsection{ Asymptotic behavior of  $cQ(m)$}

 To apply the uniform estimates obtained above to study the
asymptotic behavior of the mass functions $m(x,t),\ n(x,t)$ with
respect to the time $t$, we introduce the following lemma. The proof
is quite simple and the detail is omitted.

 \vspace{4mm}

\noindent{\bf Lemma 3.9.} \ \  Suppose that $g(t)\geq0$ for $t\geq
0$,\ $g(t)\in L^1(0, \infty)$ and $g'(t)\in L^1(0, \infty)$.\  Then
$\lim\limits_{t\to\infty}g(t)=0$.

\vspace{3mm}

Now we prove Theorem 2.2. Let
$$
g(t)=\int_0^1(cQ(m))^{\gamma}(x,t)dx. \eqno(3.29)
$$
Integrating (3.9) with respect to $x$ over $[0, 1]$, we have
$$
\int_0^{t}\int_0^1(cQ(m))^{\gamma}(x,s)dxds \leq C, \eqno(3.30)
$$
which implies $g(t)\in L^1(0, \infty)$.

 Now we prove $g'(t)\in
L^1(0, \infty)$. By the second equation of (2.7) and using
Cauchy-Schwarz inequality, we obtain
$$
\arraycolsep=1.5pt
 \begin{array}[b]{rl}
 \displaystyle\int_{0}^{\infty}|g'(t)|dt=
 & \displaystyle\gamma\int_{0}^{\infty}\left|\int_{0}^{1}(cQ(m))^{\gamma-1}c(x)Q(m)_{t}dy\right|dt \\[5mm]
 =& \displaystyle\int_{0}^{\infty}\left|\int_{0}^{1}\gamma\rho_{l}c^{\gamma}Q(m)^{\gamma+1}u_xdx\right|dt\\ [5mm]
 \leq & \displaystyle C\int_0^\infty\int_0^1c^{\beta }Q(m)^{1+\beta}u_x^2dxdt
+C\int_0^\infty\int_0^1c^{2\gamma-\beta}Q(m)^{2\gamma+1-\beta}dxdt.
\end{array}
\eqno(3.31)
$$
By $(3.8)$, $(3.9)$, and the assumptions $(A_{1})$ or $(A_{1})'$, we
can estimate the last term on the right-hand side in $(3.31)$ as
follows:
$$
\D \int_0^t\int_0^1 c^{2\gamma-\beta}Q(m)^{2\gamma+1-\beta}dxds
 \leq  \int_{0}^{t}\max \limits_{[0,1]}(cQ(m))^{\gamma}\int_{0}^{1}c^{\gamma-\beta}Q(m)^{\gamma+1-\beta}dxds
 \leq  C.$$
Substituting the above inequality into $(3.31)$ and using Lemma 3.2,
we deduce $g'(t)\in L^1(0, \infty)$.

Consequently,
$$
\lim_{t \rightarrow \infty}g(t)=0. \eqno(3.32)
$$
By (3.32), Lemma 3.3 or H$\ddot{{\rm o}}$lder inequality, we have
$$
\lim_{t \rightarrow\infty}\int_0^1(cQ(m))^{\lambda}(x, t)dx=0,
\eqno(3.33)
$$
for any $0<\lambda<\infty$.

Now we prove Theorem 2.2, which is divided into two cases:

 Case 1.
When the initial masses connect to vacuum continuously
(corresponding to the boundary condition $(2.10)$), choosing
$k>\beta>0$ and applying $(3.33)$, Lemma 3.8 and H$\ddot{{\rm
o}}$lder inequality, we have
$$\arraycolsep=1.5pt
 \begin{array}[b]{rcl}
 0\leq (cQ(m))^{k}&=&\D \int_{0}^{x}((cQ(m))^{k})_{y}dy\\[3mm]
 &=&\D \int_{0}^{x}k(cQ(m))^{k-\beta}(cQ(m))^{\beta-1}(cQ(m))_{y}dy\\[3mm]
 &=&\D \frac{k}{\beta}\int_{0}^{x}(cQ(m))^{k-\beta}((cQ(m))^{\beta})_{y}dy\\[3mm]
 &\leq &\D C \left(\int_{0}^{1}(cQ(m))^{2k-2\beta}dx\right)^{\frac{1}{2}}
 \left(\int_{0}^{1}((cQ(m))^{\beta})^{2}_{x}dx\right)^{\frac{1}{2}}\\[3mm]
 &\leq &C \D \left(\int_{0}^{1}(cQ(m))^{2k-2\beta}dx\right)^{\frac{1}{2}} \rightarrow
 0,\ \ \textrm{as}\ \ \  t\rightarrow \infty.
 \end{array}
 $$

Case 2. When the initial masses connect to vacuum discontinuously
(corresponding to the boundary condition $(2.8)$), choosing
$k>\beta>0$ and applying $(3.33)$, Corollary 3.6, Lemma 3.8 and
H$\ddot{{\rm o}}$lder inequality, we have
$$
\arraycolsep=1.5pt
 \begin{array}[b]{rcl}
 0\leq (cQ(m))^{k}&=&\D (cQ(m))^{k}(0,t)+\int_{0}^{x}((cQ(m))^{k})_{y}dy\\[3mm]
 &\leq &\D C(1+t)^{-\frac{k}{\gamma-\beta}} +\left(\int_{0}^{1}(cQ(m))^{2k-2\beta}dx\right)^{\frac{1}{2}}
 \left(\int_{0}^{1}((cQ(m))^{\beta})^{2}_{x}dx\right)^{\frac{1}{2}}\\[3mm]
 &\leq &C(1+t)^{-\frac{k}{\gamma-\beta}}+C  \D \left(\int_{0}^{1}(cQ(m))^{2k-2\beta}dx\right)^{\frac{1}{2}} \rightarrow
 0,\ \ \textrm{as}\  t \rightarrow\infty.
 \end{array}
 $$
Combining the above two cases, we have
$$
(cQ(m))(x,t)\rightarrow 0,\ \  \ \ {\rm as}\ \ \ \ t \rightarrow
\infty,
$$
which implies
$$
\lim\limits _{t\rightarrow \infty}\frac{n}{m}\cdot\frac{m}{\rho_{l}-m}=
\lim\limits _{t\rightarrow \infty}\frac{n}{\rho_{l}-m}=0.
$$
Thus
$$
\lim\limits _{t\rightarrow \infty}n(x,t)=
\lim\limits _{t\rightarrow
\infty}\frac{n}{\rho_{l}-m}\cdot(\rho_{l}-m)=0,
$$
and
$$
\lim\limits _{t\rightarrow \infty}m(x,t)
=0,
$$
for any $ x \in [0,1]$.

This completes the proof of Theorem $2.2$.

\section{Decay rates  of the mass functions }
 Now we are in the  position to estimate the stabilization  rates
of the mass functions $m(x,t),\ n(x,t)$ as $t \rightarrow \infty$.

To do this, introduce a new function $w(x,t)$ defined as follows
(cf. \cite{Nagasawa}):
$$
w(x,t)=\rho_{l}u(x, t)-\frac{1}{1+t}\int_0^x\frac{1}{Q(m)}dy
+\frac{1}{1+t}\int_0^1\int_0^x\frac{1}{Q(m)} dydx. \eqno(4.1)
$$
By direct calculation, we have
$$
w_x=\rho_{l}u_x-\frac{1}{(1+t)Q(m)}, \eqno(4.2)
$$
and
$$
w_{t}+\frac{w}{1+t}=\rho_{l}u_{t}. \eqno(4.3)
$$
Here we have used the fact that
$$
\int_0^1u(x, t)dx=\int_0^1u_0(x)dx,
$$
(see (3.1)). Assume now $\int_0^1u_0(y)dy=0$ for the simplicity of
presentation.

Thus the auxiliary functions $w$ and $Q(m)$  satisfy the following
$$
\left \{
\begin{array}{l}
\partial_{t}c=0,\\[3mm]
\D\partial_{t}Q(m)+Q(m)^{2}\partial_{x}w+\frac{Q(m)}{1+t} = 0 , \\[3mm]
\D \partial_{t}w+\frac{w}{1+t}=
\partial_{x}\left(c^{\beta}Q(m)^{\beta+1}\partial_{x}w+\frac{(cQ(m))^{\beta}}{1+t}-\rho_{l}(cQ(m))^{\gamma}\right).
\end{array}
\right. \eqno(4.4)
$$
Then we have

\vspace{4mm}

 \noindent{\bf Lemma 4.1.} \ \ Let $(c(x),u(x,
t),Q(m)(x,t))$ be a global weak solution to the fixed boundary value
problem $(2.7),\ (2.8),\ (2.9)$, or $(2.7),\ (2.10),\ (2.11)$. Then
for any $\beta>0,\ \gamma \geq 1+\beta$, the following estimates
hold for any $t>0$,

\vspace{3mm}

{\bf Case I: ${\mathbf {0<\beta<1}}$.}
$$
\arraycolsep=1.5pt
\begin{array}[b]{rl}
& \displaystyle\frac{1}{2}(1+t)^\theta \int_0^1w^2dx
+\frac{(1+t)^{\theta-1}}{1-\beta}\int_0^1c^{\beta}Q(m)^{\beta-1}dx+
\frac{\rho_{l}(1+t)^\theta}{\gamma-1}\int_0^1c^{\gamma}Q(m)^{\gamma-1}dx\\[5mm]
&\displaystyle
+\left(1-\frac{\theta}{2}\right)\int_{0}^{t}(1+s)^{\theta-1}\int_0^1w^2dxds
+\int_{0}^{t}(1+s)^\theta\int_0^1c^{\beta}Q(m)^{1+\beta}w_x^2dxds\\[5mm]
& \displaystyle
+\frac{\beta-\theta}{1-\beta}\int_{0}^{t}(1+s)^{\theta-2}\int_0^1c^{\beta}Q(m)^{\beta-1}dxds+
\rho_{l}\frac{\gamma-1-\theta}{\gamma-1}\int_{0}^{t}(1+s)^{\theta-1}\int_0^1c^{\gamma}Q(m)^{\gamma-1}dxds\\[5mm]
\leq &C,
\end{array}
\eqno(4.5)
$$
where
$$
\theta=\min\{ \gamma-1,\beta\}=\beta.\eqno(4.6)
$$

\vspace{2mm}

{\bf Case II: ${\mathbf {\beta=1}}$.}
$$
\arraycolsep=1.5pt
\begin{array}[b]{rl}
& \displaystyle
\frac{1}{2}(1+t)^\theta \int_0^1w^2dx+
\frac{\rho_{l}(1+t)^\theta}{\gamma-1}\int_0^1c^{\gamma}Q(m)^{\gamma-1}dx \\[5mm]
& \displaystyle
+\left(1-\frac{\theta}{2}\right)\int_{0}^{t}(1+s)^{\theta-1}\int_0^1w^2dxds
+\int_{0}^{1}(1+s)^\theta\int_0^1cQ(m)^{2}w_x^2dxds\\ [5mm] &
\displaystyle
+\rho_{l}\frac{\gamma-1-\theta}{\gamma-1}\int_{0}^{t}(1+s)^{\theta-1}\int_0^1c^{\gamma}Q(m)^{\gamma-1}dxds
\\[5mm]
 \leq& C+C \ln(1+t),
\end{array}
\eqno(4.7)
$$
where $\theta=1$.

\vspace{2mm}

 {\bf Case III: ${\mathbf {\beta>1}}$.}
$$
\arraycolsep=1.5pt
\begin{array}[b]{rl}
& \displaystyle
\frac{1}{2}(1+t)^\theta \int_0^1w^2dx+
\frac{\rho_{l}(1+t)^\theta}{2(\gamma-1)}\int_0^1c^{\gamma}Q(m)^{\gamma-1}dx \\[5mm]
& \displaystyle
+\left(1-\frac{\theta}{2}\right)\int_{0}^{t}(1+s)^{\theta-1}\int_0^1w^2dxds
+\int_{0}^{t}(1+s)^{\theta}\int_0^1c^{\beta}Q(m)^{1+\beta}w_x^2dxds\\ [5mm]
& \displaystyle
+\frac{\rho_{l}(\gamma-1-\theta)}{2(\gamma-1)}\int_{0}^{t}(1+s)^{\theta-1}\int_0^1c^{\gamma}Q(m)^{\gamma-1}dxds\\[5mm]
\leq & C+C(\ln(1+t))^{l},
\end{array}
\eqno(4.8)
$$
where
$$
\theta=\left\{\begin{array}{l}2,\ \ \ \ \ \ \ \ \textrm{for}\ \ \
\frac{\gamma-1}{\gamma-\beta}>2,\\[3mm]
\frac{\gamma-1}{\gamma-\beta},\ \ \ \ \ \textrm{for} \ \ \
\frac{\gamma-1}{\gamma-\beta}\leq 2,
\end{array}
\right. \eqno(4.9)
$$
and $l=0$, when $\frac{\gamma-1}{\gamma-\beta}
>2 $, whereas $l=1$, when
 $\frac{\gamma-1}{\gamma-\beta}\leq2$.

\vspace{3mm}

{\it Proof.}\ \ Multiplying $(4.4)_3$ by $w$, integrating the
resulting equation with respect to $x$ over $[0,1]$, using
integration by parts, we obtain by the boundary conditions (2.8) or
(2.10)
$$
\arraycolsep=1.5pt
\begin{array}[b]{rl}
& \displaystyle
\frac{1}{2}\frac{d}{dt}\int_0^1w^2dx+\frac{1}{1+t}\int_0^1w^2dx \\[5mm]
 = & \displaystyle \int_0^1(c^{\beta}Q(m)^{\beta+1}w_{x})_{x}wdx+
 \frac{1}{1+t}\int_0^1({c^{\beta}Q(m)^{\beta}})_{x}wdx
 -\rho_{l}\int_0^1(c^{\gamma}Q(m)^{\gamma})_xwdx\\[5mm]
 = &\displaystyle \left. c^{\beta}Q(m)^{\beta+1}w_{x}w\right|_{0}^{1}
 +\frac{1}{1+t} \left.{c^{\beta}Q(m)^{\beta}}w\right|_{0}^{1}
 -\rho_{l}c^{\gamma}Q(m)^{\gamma}w\Big|_{0}^{1}\\[5mm]
 &\displaystyle
 -\int_0^1c^{\beta}Q(m)^{\beta+1}w_{x}^{2}dx-
 \frac{1}{1+t}\int_0^1c^{\beta}Q(m)^{\beta}w_{x}dx+\rho_{l}\int_0^1c^{\gamma}Q(m)^{\gamma}w_{x}dx\\[5mm]
 =&\displaystyle
 -\int_0^1c^{\beta}Q(m)^{\beta+1}w_{x}^{2}dx-
 \frac{1}{1+t}\int_0^1c^{\beta}Q(m)^{\beta}w_{x}dx+\rho_{l}\int_0^1c^{\gamma}Q(m)^{\gamma}w_{x}dx,
\end{array}
\eqno(4.10)
$$
i.e.,
$$
\arraycolsep=1.5pt
\begin{array}[b]{rl}
& \displaystyle
\frac{1}{2}\frac{d}{dt}\int_0^1w^2dx+\frac{1}{1+t}\int_0^1w^2dx
+\int_0^1c^{\beta}Q(m)^{\beta+1}w_{x}^{2}dx\\[3mm]
=&
\displaystyle-\frac{1}{1+t}\int_0^1c^{\beta}Q(m)^{\beta}w_{x}dx+\rho_{l}\int_0^1c^{\gamma}Q(m)^{\gamma}w_{x}dx.
\end{array}
\eqno(4.11)
$$
Now we will prove (4.5), (4.7) and (4.8).

\vspace{2mm}

 {\bf Case I: ${\mathbf {0<\beta<1}}$ (The proof of
(4.5))}.

Notes that
$$
w_x=\rho_{l}u_x-\frac{1}{(1+t)Q(m)}=\left(\frac{1}{Q(m)}\right)_{t}-\frac{1}{(1+t)Q(m)}.$$
Thus we can estimate the first and second terms on the right-hand
side in (4.11) as following:
 $$
 \arraycolsep=1.5pt
\begin{array}[b]{rl}
& \displaystyle-\frac{1}{1+t}\int_0^1c^{\beta}Q(m)^{\beta}w_{x}dx\\[5mm]
 =&\displaystyle-
 \frac{1}{1+t}\int_0^1c^{\beta}Q(m)^{\beta}\left\{\left(\frac{1}{Q(m)}\right)_{t}-\frac{1}{(1+t)Q(m)}\right\}dx\\[5mm]
 =&\displaystyle-\frac{1}{(1-\beta)(1+t)}\int_0^1c^{\beta}(Q(m)^{\beta-1})_{t}dx
 +\frac{1}{(1+t)^{2}}\int_{0}^{1}c^{\beta}Q(m)^{\beta-1}dx,
\end{array}
\eqno(4.12)
$$
and
$$
\arraycolsep=1.5pt
\begin{array}[b]{rcl}
 \displaystyle\rho_{l}\int_0^1c^{\gamma}Q(m)^{\gamma}w_{x}dx&=
 &\displaystyle\rho_{l}
 \int_0^1c^{\gamma}Q(m)^{\gamma}\left\{\left(\frac{1}{Q(m)}\right)_{t}-\frac{1}{(1+t)Q(m)}\right\}dx\\[5mm]
 &=&\displaystyle\frac{\rho_{l}}{1-\gamma}\int_0^1c^{\gamma}(Q(m)^{\gamma-1})_{t}dx
 -\frac{\rho_{l}}{1+t}\int_{0}^{1}c^{\gamma}Q(m)^{\gamma-1}dx.
\end{array}
\eqno(4.13)
$$
Substituting (4.12) and (4.13) into (4.11), we get
$$
\arraycolsep=1.5pt
\begin{array}[b]{rl}
& \displaystyle
\frac{d}{dt}\int_0^1\left(\frac{w^2}{2}+
\frac{\rho_{l}}{\gamma-1}c^{\gamma}Q(m)^{\gamma-1}\right)dx
+\frac{1}{1+t}\int_0^1w^2dx\\[3mm]
&\displaystyle +\int_0^1c^{\beta}Q(m)^{\beta+1}w_{x}^{2}dx+\frac{\rho_{l}}{1+t}\int_0^1c^{\gamma}Q(m)^{\gamma-1}dx\\[3mm]
=& \displaystyle\frac{1}{(\beta-1)(1+t)}\int_0^1c^{\beta}(Q(m)^{\beta-1})_{t}dx
 +\frac{1}{(1+t)^{2}}\int_{0}^{1}c^{\beta}Q(m)^{\beta-1}dx.
\end{array}
\eqno(4.14)
$$
Multiplying (4.14) by $(1+t)^{\theta}$ for some $\theta$ to be
determined later, we deduce for any $0< \beta <1$
$$
\arraycolsep=1.5pt
\begin{array}[b]{rl}
& \displaystyle
\frac{d}{dt}\left\{\frac{1}{2}(1+t)^\theta \int_0^1w^2dx
+\frac{(1+t)^{\theta-1}}{1-\beta}\int_0^1c^{\beta}Q(m)^{\beta-1}dx+
\frac{\rho_{l}(1+t)^\theta}{\gamma-1}\int_0^1c^{\gamma}Q(m)^{\gamma-1}dx\right\} \\
[5mm] & \displaystyle
+\left(1-\frac{\theta}{2}\right)(1+t)^{\theta-1}\int_0^1w^2dx
+(1+t)^\theta\int_0^1c^{\beta}Q(m)^{1+\beta}w_x^2dx
\\ [5mm]
& \displaystyle
+\frac{\beta-\theta}{1-\beta}(1+t)^{\theta-2}\int_0^1c^{\beta}Q(m)^{\beta-1}dx+
\rho_{l}\frac{\gamma-1-\theta}{\gamma-1}(1+t)^{\theta-1}\int_0^1c^{\gamma}Q(m)^{\gamma-1}dx
\\[5mm]
 = & 0.
\end{array}
\eqno(4.15)
$$
Taking $\theta=\min \ \{\beta, \gamma-1\}=\beta$ in (4.15) and
integrating (4.15) with respect to $t$ over $[0, t]$, we deduce
(4.5).

Consequently,
$$
\int_{0}^{1}c^{\gamma}Q(m)^{\gamma-1}dx\leq
C(1+t)^{-\theta}.\eqno(4.16)
$$

{\bf Case II: ${\mathbf{\beta=1}}$ (The proof of (4.7))}. Under this
case, the first term on the right-hand side in (4.11) can be
rewritten as
$$
\arraycolsep=1.5pt
\begin{array}[b]{rl}
&  \displaystyle-\frac{1}{1+t}\int_0^1cQ(m)w_{x}dx\\[5mm]
=& \displaystyle
-\frac{1}{1+t}\int_{0}^{1}cQ(m)\left\{\left(\frac{1}{Q(m)}\right)_{t}-\frac{1}{(1+t)Q(m)}\right\}dx
\\ [5mm] = & \displaystyle
\frac{1}{(1+t)}\int_0^1c(\ln(Q(m))_{t}dx+\frac{1}{(1+t)^2}\int_{0}^{1}c_{0}(x)dx.
\end{array}
\eqno(4.17)
$$
Similar to (4.15), we have:
$$
\arraycolsep=1.5pt
\begin{array}[b]{rl}
& \displaystyle
\frac{d}{dt}\left\{\frac{1}{2}(1+t)^\theta \int_0^1w^2dx+
\frac{\rho_{l}(1+t)^\theta}{\gamma-1}\int_0^1c^{\gamma}Q(m)^{\gamma-1}dx\right\} \\[5mm]
& \displaystyle
+\left(1-\frac{\theta}{2}\right)(1+t)^{\theta-1}\int_0^1w^2dx
+(1+t)^\theta\int_0^1cQ(m)^{2}w_x^2dx\\ [5mm] & \displaystyle
+\rho_{l}\frac{\gamma-1-\theta}{\gamma-1}(1+t)^{\theta-1}\int_0^1c^{\gamma}Q(m)^{\gamma-1}dx
\\[5mm]
 =&\displaystyle \frac{d}{dt}\left\{(1+t)^{\theta-1}\int_{0}^{1}c\ln(Q(m))dx  \right\}
 +(1+t)^{\theta-2}\int_{0}^{1}c_{0}(x)dx\\[5mm]
 &\displaystyle+(1-\theta)(1+t)^{\theta-2}\int_{0}^{1}c\ln(Q(m))dx.
\end{array}
\eqno(4.18)
$$
By using $\ln x\leq x-1$ for any $x>0$ and Lemma 3.3, we have
$$\int_{0}^{1}c\ln Q(m)dx \leq \int_{0}^{1}cQ(m)dx\leq C$$
and the assumption $(A_1)$ or $(A_1)'$ implies that $$\int
_{0}^{1}c_{0}(x)dx \leq C.$$ Taking $\theta=1$ in (4.18) and
integrating (4.18) with respect to $t$ over $[0, t]$, we deduce
(4.7).

Consequently,
$$
\int_{0}^{1}c^{\gamma}Q(m)^{\gamma-1}dx\leq
C(1+t)^{-\theta}\ln(1+t).\eqno(4.19)
$$

\vspace{2mm}

{\bf Case III: ${\mathbf{\beta>1}}$ (The proof of (4.8))}.

Rewrite (4.15) as
$$
\arraycolsep=1.5pt
\begin{array}[b]{rl}
& \displaystyle
\frac{d}{dt}\left\{\frac{1}{2}(1+t)^\theta \int_0^1w^2dx+
\frac{\rho_{l}(1+t)^\theta}{\gamma-1}\int_0^1c^{\gamma}Q(m)^{\gamma-1}dx\right\} \\[5mm]
& \displaystyle
+\left(1-\frac{\theta}{2}\right)(1+t)^{\theta-1}\int_0^1w^2dx
+(1+t)^\theta\int_0^1c^{\beta}Q(m)^{1+\beta}w_x^2dx\\ [5mm]
& \displaystyle
+\rho_{l}\frac{\gamma-1-\theta}{\gamma-1}(1+t)^{\theta-1}\int_0^1c^{\gamma}Q(m)^{\gamma-1}dx
\\[5mm]
 =&\displaystyle \frac{d}{dt}\left\{\frac{(1+t)^{\theta-1}}{\beta-1}\int_{0}^{1}c^{\beta}Q(m)^{\beta-1}dx \right\}
 +\frac{\beta-\theta}{\beta-1}(1+t)^{\theta-2}\int_{0}^{1}c^{\beta}Q(m)^{\beta-1}dx.
\end{array}
\eqno(4.20)
$$
Integrating (4.20) with respect to $t$ over $[0,t]$, we have
$$
\arraycolsep=1.5pt
\begin{array}[b]{rl}
& \displaystyle
\frac{1}{2}(1+t)^\theta \int_0^1w^2dx+
\frac{\rho_{l}(1+t)^\theta}{\gamma-1}\int_0^1c^{\gamma}Q(m)^{\gamma-1}dx \\[5mm]
& \displaystyle
+\left(1-\frac{\theta}{2}\right)\int_{0}^{t}(1+s)^{\theta-1}\int_0^1w^2dxds
+\int_{0}^{t}(1+s)^{\theta}\int_0^1c^{\beta}Q(m)^{1+\beta}w_x^2dxds\\
[5mm] & \displaystyle
+\rho_{l}\frac{\gamma-1-\theta}{\gamma-1}\int_{0}^{t}(1+s)^{\theta-1}\int_0^1c^{\gamma}Q(m)^{\gamma-1}dxds
\\[5mm]
=&\displaystyle \frac{1}{2}\int_0^1w_{0}^2dx+
\frac{\rho_{l}}{\gamma-1}\int_0^1c_{0}^{\gamma}Q(m_{0})^{\gamma-1}dx-\frac{1}{\beta-1}\int_{0}^{1}c_{0}^{\beta}Q(m_{0})^{\beta-1}dx  \\[5mm]
&+\displaystyle
\frac{(1+t)^{\theta-1}}{\beta-1}\int_{0}^{1}c^{\beta}Q(m)^{\beta-1}dx
 +\frac{\beta-\theta}{\beta-1} \int_{0}^{t}(1+s)^{\theta-2}\int_{0}^{1}c^{\beta}Q(m)^{\beta-1}dxds\\[5mm]
=&\displaystyle-\frac{1}{\beta-1}\int_{0}^{1}c_{0}^{\beta}Q(m_{0})^{\beta-1}dx
+\frac{1}{2}\int_0^1w_{0}^2dx+
\frac{\rho_{l}}{\gamma-1}\int_0^1c_{0}^{\gamma}Q(m_{0})^{\gamma-1}dx+I_{1}+I_{2} . \\[5mm]
\end{array}
\eqno(4.21)
$$
By Young inequality, we have
$$
\arraycolsep=1.5pt
\begin{array}[b]{rcl}
I_{1}&=&\displaystyle \frac{(1+t)^{\theta-1}}{\beta-1}\int_{0}^{1}c^{\beta}Q(m)^{\beta-1}dx\\[5mm]
&=&\displaystyle
\frac{1}{\beta-1}\int_{0}^{1}(cQ(m))^{\beta-1}(1+t)^{\frac{(\beta-1)\theta}{\gamma-1}}c
(1+t)^{\theta-1-\frac{(\beta-1)\theta}{\gamma-1}}dx\\[5mm]
& \leq &\displaystyle
\frac{\rho_{l}(1+t)^{\theta}}{2(\gamma-1)}\int_{0}^{1}c^{\gamma}Q(m)^{\gamma-1}dx
+C\displaystyle (1+t)^{(\theta-1-\frac{(\beta-1)\theta}{\gamma-1})\frac{\gamma-1}{\gamma-\beta}}\int_{0}^{1}c_{0}(x)dx\\[5mm]
& \leq &\displaystyle
\frac{\rho_{l}(1+t)^{\theta}}{2(\gamma-1)}\int_{0}^{1}c^{\gamma}Q(m)^{\gamma-1}dx
+C\displaystyle
(1+t)^{(\theta-1-\frac{(\beta-1)\theta}{\gamma-1})\frac{\gamma-1}{\gamma-\beta}},
\end{array}
\eqno(4.22)
$$
and
$$
\arraycolsep=1.5pt
\begin{array}[b]{rcl}
I_{2}&=&\displaystyle \frac{\beta-\theta}{\beta-1}\int_{0}^{t}(1+s)^{\theta-2}\int_{0}^{1}c^{\beta}Q(m)^{\beta-1}dxds\\[5mm]
& \leq &\displaystyle \frac{\rho_{l}(\gamma-1-\theta)}{2(\gamma-1)}
\int_{0}^{t}(1+s)^{\theta-1}\int_{0}^{1}c^{\gamma}Q(m)^{\gamma-1}dx
+C\displaystyle\int_{0}^{t}
(1+s)^{(\theta-2-\frac{(\beta-1)(\theta-1)}{\gamma-1})\frac{\gamma-1}{\gamma-\beta}}ds.
\end{array}
\eqno(4.23)
$$
Substituting (4.22) and (4.23) into (4.21), we have
$$
\arraycolsep=1.5pt
\begin{array}[b]{rl}
& \displaystyle
\frac{1}{2}(1+t)^\theta \int_0^1w^2dx+
\frac{\rho_{l}(1+t)^\theta}{2(\gamma-1)}\int_0^1c^{\gamma}Q(m)^{\gamma-1}dx \\[5mm]
& \displaystyle
+\left(1-\frac{\theta}{2}\right)\int_{0}^{t}(1+s)^{\theta-1}\int_0^1w^2dxds
+\int_{0}^{t}(1+s)^{\theta}\int_0^1c^{\beta}Q(m)^{1+\beta}w_x^2dxds\\
[5mm] & \displaystyle
+\frac{\rho_{l}(\gamma-1-\theta)}{2(\gamma-1)}\int_{0}^{t}(1+s)^{\theta-1}\int_0^1c^{\gamma}Q(m)^{\gamma-1}dxds
\\[5mm]
\leq& C+C\displaystyle
(1+t)^{\theta-\frac{\gamma-1}{\gamma-\beta}}+C\int_{0}^{t}
(1+s)^{\theta-1-\frac{\gamma-1}{\gamma-\beta}}ds.
\end{array}
\eqno(4.24)
$$
Taking $\theta=2$, when $\frac{\gamma-1}{\gamma-\beta}>2$ in (4.24),
we have
$$
\int_{0}^{1}c^{\gamma}Q(m)^{\gamma-1}dx\leq
C(1+t)^{-\theta}.\eqno(4.25)$$ Taking
$\theta=\frac{\gamma-1}{\gamma-\beta}$, when
$\frac{\gamma-1}{\gamma-\beta}\leq2 $ in (4.24), we have
$$
\int_{0}^{1}c^{\gamma}Q(m)^{\gamma-1}dx\leq
C(1+t)^{-\theta}\ln(1+t).\eqno(4.26)
$$
This completes the proof of
Lemma 4.1.

 \vspace{4mm}

\textbf{The proof of Theorem 2.3}. \ \ Under the boundary condition
$(2.8)$, for $0< \beta <1$ or $\beta>1, \
\frac{\gamma-1}{\gamma-\beta}>2$, choosing some constant
$2k=\gamma-1+2\beta$ and using the assumption $(A_{1})$, Corollary
3.6, Lemma 3.8, $(4.16)$ and $(4.25)$, we have
$$
\arraycolsep=1.5pt
\begin{array}[b]{rl}
\displaystyle (cQ(m))^k(x, t)= & \displaystyle \ (cQ(m))^k(0,t)
+\int_0^x((cQ(m))^k)_{y}(y, t)dy\\ [3mm]
 \leq & \displaystyle C(1+t)^{-\frac{k}{\gamma-\beta}}
+C \left(\int_0^1((cQ(m))^{\beta})_{x}^{2}dx\right)^{\frac{1}{2}}
\left(\int_0^1(cQ(m))^{2k-2\beta}dx\right)^{\frac{1}{2}} \\ [3mm]
\leq&\displaystyle C(1+t)^{-\frac{k}{\gamma-\beta}}+
C\left(\int_0^1\frac{1}{c_{0}}  c^{\gamma}Q(m)^{\gamma-1}dx\right)^{\frac{1}{2}}\\
[3mm]
 \leq& \displaystyle
C(1+t)^{-\frac{k}{\gamma-\beta}}+C(1+t)^{-\frac{\theta}{2}}\\ [3mm]
\leq &\displaystyle C(1+t)^{-\frac{\theta}{2}},
\end{array}
\eqno(4.27)
$$
which  implies
$$
(cQ(m))(x,t)\leq C(1+t)^{-\frac{\theta}{2k}}=
C(1+t)^{-\frac{\theta}{\gamma-1+2\beta}}, \eqno(4.28)
$$
i.e.,
$$
\frac{ n(x,t)}{\rho_{l}-m(x,t)}\leq
C(1+t)^{-\frac{\theta}{\gamma-1+2\beta}}. \eqno(4.29)
$$
Thus
$$
 n(x,t)=\frac{ n(x,t)}{\rho_{l}-m(x,t)} \cdot (\rho_{l}-m(x,t))
 \leq C(1+t)^{-\frac{\theta}{\gamma-1+2\beta}}, \eqno(4.30)
$$
and
$$
 m(x,t)= n(x,t) \cdot c(x)^{-1}
 \leq C(1+t)^{-\frac{\theta}{\gamma-1+2\beta}}, \eqno(4.31)
$$
for  any $ x \in [0,1]$.

Similarly, if $\beta=1$ or $\beta>1,\ \frac{\gamma-1}{\gamma-\beta}
\leq2$, we have
$$
(cQ(m))^k(x, t)\leq \displaystyle C(1+t)^{-\frac{\theta}{2}}\sqrt{\ln(1+t)},
$$
which implies
$$
n(x,t)\leq
C(1+t)^{-\frac{\theta}{\gamma-1+2\beta}}(\ln(1+t))^{\frac{1}{\gamma-1+2\beta}},
\eqno(4.32)
$$
and
$$
m(x,t)\leq
C(1+t)^{-\frac{\theta}{\gamma-1+2\beta}}(\ln(1+t))^{\frac{1}{\gamma-1+2\beta}},
\eqno(4.33)
$$
for any $ x \in [0,1]$. Here we have used $(4.19)$ and $(4.26)$.

Under the boundary condition (2.10), for $ 0< \beta <1$ or
$\beta>1,\ \frac{\gamma-1}{\gamma-\beta}>2$, choosing some constant
$2k_{1}=\frac{\gamma-1}{2}+2\beta$ and using the assumptions
$(A_{1})'$, Lemma 3.8, (4.16), (4.25) and H\"{o}lder's inequality,
we have
$$
\arraycolsep=1.5pt
\begin{array}[b]{rl}
\displaystyle (cQ(m))^{k_{1}}(x, t)= &
\displaystyle\int_0^x((cQ(m))^{k_{1}})_{y}(y, t)dy\\ [5mm]
 \leq & \displaystyle C \left(\int_0^1((cQ(m))^{\beta})_{x}^{2}dx\right)^{\frac{1}{2}}
\left(\int_0^1(cQ(m))^{2k_{1}-2\beta}dx\right)^{\frac{1}{2}} \\
[5mm] \leq& C
\displaystyle\left(\int_0^1\frac{1}{c^{\frac{1}{2}}}\cdot
c^{\frac{\gamma}{2}}Q(m)^{\frac{\gamma-1}{2}}dx\right)^{\frac{1}{2}}\\[5mm]
\leq& C \displaystyle \left(\int_0^1
c^{\gamma}Q(m)^{\gamma-1}dx\right)^{\frac{1}{4}}
\left(\int_0^1\frac{1}{ c_{0}(x)}dx\right)^{\frac{1}{4}}\\[5mm]
 \leq& \displaystyle C(1+t)^{-\frac{\theta}{4}},\\ [5mm]
\end{array}
\eqno(4.35)
$$
which implies
$$ (cQ(m))(x,t)\leq C(1+t)^{-\frac{\theta}{4k_{1}}}=
C(1+t)^{-\frac{\theta}{\gamma-1+4\beta}}, \eqno(4.36)
$$
i.e.,
$$
\frac{ n(x,t)}{\rho_{l}-m(x,t)}\leq
C(1+t)^{-\frac{\theta}{\gamma-1+4\beta}}, \eqno(4.37)
$$
Thus
$$
 n(x,t)=\frac{ n(x,t)}{\rho_{l}-m(x,t)} \cdot (\rho_{l}-m(x,t))
 \leq C(1+t)^{-\frac{\theta}{\gamma-1+4\beta}}, \eqno(4.38)
$$
and
$$
 m(x,t)
 \leq C(1+t)^{-\frac{\theta}{\gamma-1+4\beta}}, \eqno(4.39)
$$
for any  $x \in [0,1]$.

 Similarly, if $\beta=1$ or $\beta>1,\
\frac{\gamma-1}{\gamma-\beta} \leq 2$, we have
$$(cQ(m))^{k_{1}}(x, t)\leq \displaystyle C(1+t)^{-\frac{\theta}{4}}(\ln(1+t))^{\frac{1}{4}},$$
then
$$
 n(x,t)\leq C(1+t)^{-\frac{\theta}{\gamma-1+4\beta}}(\ln(1+t))^{\frac{1}{\gamma-1+4\beta}}, \eqno(4.40)
$$
and
$$
 m(x,t)\leq C(1+t)^{-\frac{\theta}{\gamma-1+4\beta}}(\ln(1+t))^{\frac{1}{\gamma-1+4\beta}}, \eqno(4.41)
$$
for any $ x \in [0,1]$. Here we have used $(4.19)$ and $(4.26)$.

 The proof of Theorem 2.3 is completed.

 \vskip 1.2cm

{\bf Acknowledgement:} \ \ The research was supported by the
National Natural Science Foundation of China $\#$10625105,
$\#$11071093, the PhD specialized grant of the Ministry of Education
of China $\#$20100144110001, and the self-determined research funds
of CCNU from the colleges'basic research and operation of MOE.

\bibliographystyle{plain}

\end{document}